\documentclass{article}

\usepackage{epsfig}
\usepackage{graphics}

\usepackage{amsmath}
\usepackage{amsfonts}
\usepackage{amssymb}

\numberwithin{equation}{section}

\newcommand{\D}[1]{{\mathbb#1}}% Doubled -Blackboard bold - caps only

\newcommand{\RR}{{\D{R}}}

\def\om{\omega}
\def\a{\alpha}

\def\lV{\left\Vert }
\def\rV{\right\Vert }

\def\Ga{\Gamma}

\def\f{\mathcal F}
\def\u{\mathcal U}
\def\proj{\operatorname{proj}}

\newtheorem{defn}{Definition}[section]

\newtheorem{thm}[defn]{Theorem}

\newtheorem{cor}[defn]{Corollary}
\newtheorem{lem}[defn]{Lemma}
\newtheorem{prop}[defn]{Proposition}

\newtheorem{rem}[defn]{Remark}

\begin{document}
\bibliographystyle{plain}

\title{Liao Standard Systems and  Nonzero Lyapunov Exponents for
Differential Flows}

\author{Wenxiang Sun\footnote{ Supported partly by NNSFC 10171004.
The first author also thanks Ohio University for its hospitality
during the winter and spring of 2002 when this paper was written.}\\
School of Mathematical Sciences\\
Peking University\\
Beijing 100871, China\\
sunwx@math.pku.edu.cn\\
\and Todd Young\\
Department of Mathematics\\
Ohio University\\
Athens 45701, Ohio, USA\\
young@math.ohiou.edu\\
fax: 740-593-9805}

\date{\today}

\maketitle

\centerline{This paper is dedicated to the memory of Professor Liao Shantao, 1920--1997.}

\begin{abstract}
Consider  a $C^1$ vector field together with
an ergodic invariant probability that has $\ell$ nonzero Lyapunov exponents.
Using orthonormal moving frames along certain transitive orbits we construct  a
linear system of  $\ell$ differential equations which is a reduced form of
Liao's ``standard system".  We show that the Lyapunov exponents
of this linear system coincide with all the  nonzero exponents of the
given vector field with respect to the given
probability. Moreover, we prove that these Lyapunov exponents
have a persistence property that implies that a
``Liao perturbation"  preserves both sign and value
of nonzero Lyapunov exponents.
\end{abstract}


\medskip


\medskip

\noindent
{\bf Running title:} Liao Systems and Lyapunov Exponents.

\newpage

\section{Introduction}

Lyapunov exponents measure the asymptotic exponential rate at which
infinitesimally nearby points approach or move away from each other
as time increases to infinity. For  a  uniformly  hyperbolic system
with positive (resp. negative) Lyapunov exponents, its nearby system
has positive (resp. negative) Lyapunov exponents as well.
Using orthonormal frames moving along certain transitive orbits
Liao (see \cite{lia96}) constructed a system of linear equations,
known as a ``standard system". In the hyperbolic case, our Main Theorem
together with a result of Liao's \cite[Theorem~2.4.1]{lia96} shows
that the Lyapunov exponents of the standard system
coincide with those of the original flow. Professor Liao had
conjectured this result. For the complement of
uniform hyperbolicity in the space of all $C^1$ systems with $C^1$
topology,   understanding dynamics through  Lyapunov exponents
and SRB measures is   incomplete but very important (see Palis \cite{Pal}).
Young \cite{You93,You97}  constructed open sets of
nonuniform hyperbolicity
cocycles for  certain special systems. In \cite{Via} Viana constructed an
open set  of  systems with multidimensional nonhyperbolic attractors
which have  SRB measures \cite{alv}.
For a compact surface, Bochi \cite{boc} showed that there is a residual set of
$C^1$ area  preserving diffeomorphisms so that each diffeomorphism in
the set is either Anosov or has a  zero Lyapunov exponent almost everywhere.

In the 1960's, Liao (see \cite{lia96}) constructed a system of linear equations,
known as a ``standard system". This system is essentially the 
variational equations along a typical orbit with respect to a typical
orthonormal frame evolving along the orbit. Liao had used the
standard system to give independent proofs of the $C^1$ closing
lemma \cite[Appendix A]{lia96} and of the topological stability 
for Anosov flows [7, Chapter 2]. While Liao's approach
is obviously philosophically related to Lyapunov exponents, the
connection has never been rigorously shown.
In the hyperbolic case, our Main Theorem
together with a result of Liao's \cite[Theorem~2.4.1]{lia96} shows
that the Lyapunov exponents of the standard system
coincide with those of the original flow. Professor Liao had
conjectured this result.

We work with  $C^1$
vector fields and develop a reduced form of Liao's standard systems.
We consider a $C^1$ vector field together with an ergodic invariant
probability that has $\ell$ nonzero Lyapunov exponents.
Using typical moving orthonormal $\ell$-frames along typical transitive orbits
of the ergodic measure,
and by using a characterization of the Lyapunov spectrum \cite{lia73,Sun98a}
we construct  a ``reduced standard system" of differential equations and
show that its  Lyapunov exponents  coincide with  the  nonzero exponents
of the original vector field. In the final section we show that the nonzero
Lyapunov exponents of the reduced standard system have certain persistence
properties.

Now let us describe the main theorem of the present paper.
We denote by $M^n$ a compact smooth $n$-dimensional Riemannian manifold and
by $S$ a $C^1$ differential system, or in other words, a $C^1$ vector field on
$M^n$. As usual $S$   induces a one-parameter transformation group
$\phi_t\colon M^n\to M^n$, $t\in \RR$ on the state manifold and
therefore a one-parameter
transformation group $\Phi_t=d\phi_t\colon TM^n\to TM^n$, $t\in \RR$ on the
tangent bundle. A probability $\nu$ on $M^n$ is $\phi$-invariant if it is
$\phi_t$-invariant for any $t\in \RR$. A $\phi$-invariant probability is
called $\phi$-ergodic if every $\phi$-invariant set has zero or full probability.
For a compact metric space $X$ and a topological flow $\varphi_t$
on it we denote by $E(X, \varphi)$ the set of all $\phi$-invariant and ergodic
probabilities.  Let $\nu$ be a $\phi$-invariant and
ergodic probability, i.e., $\nu\in E(M^n, \phi)$. From the  Multiplicative
Ergodic Theorem (see \cite{lia63,Ole}), there
exists a $\phi_t$-invariant subset $B$, with $\nu$-full probability, such
that for any $x\in B$ and $u\in T_xM^n$ the following limit, called Lyapunov
exponent, exists:
\begin{equation}\label{def:lyap}
\lambda:=\lim_{t\to \infty}\frac 1t\log\lV \Phi_t(u)\rV\,\,(\text{ or }
\lambda:=\lim_{t\to -\infty}\frac 1t\log\lV \Phi_t(u)\rV ).
\end{equation}
It is known that $\nu$ has at most $n$ different Lyapunov exponents, where $n$ indicates the
dimension of the state manifold $M^n$.

\bigskip

\noindent
{\bf  Main Theorem} \ \
{\em Suppose that a $\phi$-invariant and ergodic probability,
$\nu\in E(M^n, \phi)$,
has $\ell$ simple nonzero Lyapunov exponents
\begin{equation}\label{nzspec}
\lambda_1 < \lambda_2 < \ldots < \lambda_\ell,
\end{equation}
together with $n-\ell$ zero Lyapunov exponents.
Then the reduced standard linear system (defined in Section 4),
\begin{equation}\label{rsls1}
\frac{dy}{dt}=y A_{\ell\times \ell}(t), \quad y\in \RR^\ell, \quad t\in \RR,
\end{equation}
is well defined and has the following properties:
\begin{enumerate}
\item The matrix $A_{\ell\times \ell}(t)$ is uniformly bounded
and continuous with respect to $t$.

\item There exist $u_1, u_2,  \ldots , u_\ell\in \RR^\ell$ such that
$$
\lim_{t\to \infty} \frac{1}{t} \log \bigl\|  y(t, u_i) \bigr\|  =\lambda_i,
$$
where $y(t,v)$ denotes a unique solution of the initial value problem
\begin{equation}\label{rsls:ivp}
\frac{dy}{dt} = y A_{\ell\times \ell}(t), \qquad y(0,v) = v.
\end{equation}

\item Consider a perturbation of the linear system
\begin{equation}\label{pert}
\frac{dy}{dt}=y A_{\ell\times \ell}(t)\,+f(t, y), \quad \sup_{t\in \RR, y\in \RR^\ell}\lV
f(t, y)\rV\leq L<\infty,
\end{equation}
where $f(t, y)$ is Lipschitz in $y$.
Then there exist $u_1^*, u_2^*,  \ldots , u_\ell ^*\in \RR^\ell$ such that
$$
\lim_{t \to \infty} \frac{1}{t} \log \bigl\|  y(t, u_i^*) \bigr\|  = \lambda_i,
$$
where $y(t, v)$ denotes a unique solution of the initial value  problem
\begin{equation}\label{pert:ivp}
\frac{dy}{dt}=y A_{\ell\times \ell}(t)\,+f(t, y), \qquad y(0,v)=v.
\end{equation}
\end{enumerate}
}

In Section 2 we review  frame bundles  and the corresponding
one parameter transformation groups induced by a given vector field.
In Section 3 we construct a reduced version of Liao's
``qualitative functions"  and then use them
to present a characterization of the Lyapunov spectrum. In Section
4 we construct the reduced standard  linear system  of  $\ell$  differential
equations  on a given probability and
establish a relation between the nonzero Lyapunov exponents of this
probability and
that of the linear system. We complete the proof of the Main Theorem  in
Section 5. An example  in Section 5 illustrates that
the original standard linear system of  $n$  differential equations introduced
by Liao \cite[Chapter 2]{lia96}
fails to satisfy the conclusions of the Main Theorem, and so, it is necessary
to develop the reduced standard linear system of
$\ell$ differential equations for the Main Theorem.
In Section 7 we present  the notion of Liao perturbation and point out by
the Main Theorem that a certain type of perturbation, known as ``Liao perturbation",
preserves the nonzero Lyapunov exponents.

%%%%%%%%%%%%%%%%%%%%%%%%%%%%%%%%%%%%%%%%%%%%%
%%%%%%%%%%%%%%%%%%%%%%%%%%%%%%%%%%%%%%%%%%%%%
%%%%%%%%%%%%%%%%%%%%%%%%%%%%%%%%%%%%%%%%%%%%%

\section{One parameter transformation groups}

We start from a $C^1$ vector field $S$  on  a compact smooth $n$-dimensional
Riemannian
manifold  $M^n$, and its induced  one-parameter transformation groups
$\phi_t\colon M^n\to M^n$, $t\in \RR$ on the state manifold and
$\Phi_t =d\phi_t\colon TM^n\to TM^n$, on the tangent bundle.

Fix some integer $\ell$, $1 \le \ell \le n$. Construct a bundle
$\u_\ell = \bigcup_{x\in M^n}\u_\ell(x)$ of $\ell$-frames, where the fiber over $x$ is
\begin{equation}\label{uell}
\u_\ell(x) = \{(u_1,\dots,u_\ell)\in T_xM^n\times \dots \times
T_xM^n : u_1,u_2,\dots,u_\ell, \text{ are linearly independent}\}.
\end{equation}

Let $p_\ell\colon\u_\ell\to M^n$ denote the bundle projection. Denote by
$\text{proj}_k\colon\u_\ell\to TM^n$  the map which sends $\a\in\u_\ell$ to
the $k$-th vector in $\a$. The vector field $S$ induces a one-parameter
transformation group on $\u_\ell$,
which we denote (with the same notation as the tangent map for the sake of
simplicity) by $\Phi_t$, $t\in \RR$, namely,
$$
\Phi_t(u_1,u_2,\dots,u_\ell)=(d\phi_t(u_1),\, d\phi_t(u_2),\dots,d\phi_t(u_\ell)).
$$
For $\a=(u_1,u_2,\dots,u_\ell)\in\u_\ell$ and a
nondegenerate $\ell \times\ell$ matrix $B=(b_{ij})$ we write
$$
\a\circ B=\left(\sum_{i=1}^\ell b_{i1}u_i, \sum_{i=1}^\ell b_{i2}u_i,\dots,
\sum_{i=1}^\ell a_{i\ell}u_i \right).
$$
Then $\Phi_t(\a\circ B)=\Phi_t(\a)\circ B$. By the Gram-Schmidt orthogonalization
process there exists a unique upper triangular matrix $\Ga(\a)$ with diagonal
elements $1$  such that $\a\circ\Ga(\a)$ is orthogonal.

Construct the  bundle $\f_\ell=U_{x\in M^n} \f_\ell(x)$ of $\ell$-orthogonal
frames, where the fiber over $x$ is
\begin{equation}\label{fell}
\f_\ell(x)= \{(u_1,u_2,\dots,u_\ell)\in\u_\ell(x)\mid \langle u_i,u_j\rangle=
0,\,\, 1\le i\ne j\le \ell\}.
\end{equation}
The bundle projection is given by $q_\ell=p_\ell|\f_\ell$. The vector field $S$
then induces a one-parameter transformation group
\begin{equation}\label{chit}
\chi_t\colon\f_\ell\to\f_\ell:\a\mapsto\Phi_t(\a)\circ\Ga(\Phi_t(\a)).
\end{equation}
If we define $\pi\colon\u_\ell\to\f_\ell$ by $\a\mapsto\a\circ\Ga(\a)$ then
$\chi_t(\a)=\pi(\Phi_t(\a))$.

Construct a bundle $\f_\ell^{\#}=U_{x\in M^n}\f_\ell^{\#}(x)$ of orthonormal
$\ell$-frames, where the fiber over $x$ is
\begin{equation}\label{felln}
\f_\ell^{\#}(x)=\{(u_1,u_2,\dots,u_\ell)\in\f_\ell(x)\mid \lV u_i \rV=1,\,\,
i=1,2,\dots,\ell\}.
\end{equation}
Then $\f_\ell^{\#}$ is a compact metrizable space. Let
$\pi^{\#}\colon\f_\ell\to\f_\ell^{\#}$ be given by
$$
\pi^{\#}(u_1,u_2,\dots,u_\ell)= \left( \frac{u_1}{\lV u_1 \rV},
\frac{u_2}{\lV u_2 \rV},\dots,\frac{u_\ell}{\lV u_\ell \rV} \right).
$$
Setting $\chi_t^{\#}=\pi^{\#}\circ(\chi_t|\f_\ell^{\#})$, we get a
one-parameter transformation group $\chi_t^{\#}\colon\f_\ell^{\#}\to\f_\ell^{\#}$.
Let $q_\ell^{\#}=q_\ell|\f_\ell^{\#}$, then $q_\ell^{\#}$ is a bundle
projection. It is easy to check that the following properties hold:
\begin{equation}\label{commute}
q_\ell\circ\chi_t=\phi_t\circ q_\ell, \quad
q_\ell^{\#}\circ\chi_t^{\#}=\phi_t\circ q_\ell^{\#}, \quad
\chi_t^{\#}\circ\pi^{\#}=\pi^{\#}\circ\chi_t.
\end{equation}

\begin{rem}
We point out that the unitary $\ell$-bundle $\u_\ell^{\#}$
of $\u_\ell$ is not necessarily  a compact metric space. For instance, when $\ell=2$, there are
sequences of 2-frames $\{\a^m=(u_1^m, u_2^m)\}_{m=1}^\infty$ such that the angle
between $u_1^m$ and $u_2^m$ goes to zero as $m \to +\infty$. Such a sequence of
frames has no accumulating point inside $\u_2^{\#}$.
\end{rem}

%%%%%%%%%%%%%%%%%%%%%%%%%%%%%%%%%%%%%%%%%%%%%%%%
%%%%%%% Qualitative Functions %%%%%%%%%%%%%%%%%%%%%%%%%%%%
%%%%%%%%%%%%%%%%%%%%%%%%%%%%%%%%%%%%%%%%%%%%%%%%

\section{Qualitative functions}

For $\a\in\f_n$, let $\zeta_{\a k}(t)=\lV \text{proj}_k\chi_t(\a) \rV$,
$k=1,2,\dots,n$. Note that $\zeta_{\a k}(t)>0$ for any $t\in \RR$.

\begin{defn}
For each $k=1,2, \ldots, n$, we call $\omega_k$ defined by:
\begin{equation}\label{omegak}
\omega_k\colon\f_n\to \RR: \a \mapsto \frac{d\zeta_{\a k}(t)}{dt}\bigl|_{t=0}
\end{equation}
a qualitative function over the orthogonal $n$-frame bundle $\f_n$ and call
$\om_k|\f_n^{\#}$ a qualitative function over the orthonormal $n$-frame bundle
$\f^{\#}_n$.
\end{defn}

The qualitative function for vector fields was introduced by Liao in
1963,  and it  plays an important role in Liao  theory  [4-8].
Sun introduced its diffeomorphism version
and described its relation with Lyapunov exponents in \cite{Sun98b,Sun99}, and
determined in \cite{Sun01} the entropy of certain classes of Grassmann bundle systems by using
these functions.

From the  definition it is easy to show  that $\om_k(\a)$ is continuous,
$\om_k(\chi_t(\a))= \frac{d\zeta_{\a k}(t)}{dt}$, and
$\om_k(\chi^{\#}_t(\a))=\frac1{\zeta_{\a k}(t)} \frac{d\zeta_{\a k}(t)}
{dt}$, so the following lemma is clear.
\begin{lem}\label{logint}
For $\a\in\f_n^{\#}$ and $k=1,2,\dots,n$, we have that $\log\zeta_{\a k}(T)=\int_0^T\om_k(\chi^{\#}_t(\a))\, dt$.
\end{lem}

If we denote by $Q_\nu(M^n, \phi)$ the set of all points $x\in M^n$ that
satisfy,
for any continuous function $f$ on $M^n$,
\begin{equation}\label{avg}
\lim_{t\to \pm\infty}\frac 1t\int_0^tf(\phi_\tau(x))d\tau =\int_{M^{n}}fd\mu,
\end{equation}
then  $Q_\nu(M^n, \phi)$ is $\phi$-invariant subset with $\nu$-full
probability. Similarly one can define
$Q_\mu(\f^{\#}_n ,  \chi^{\#} )  $ for any probability $\mu\in
E(\f^{\#}_n, \chi^{\#})$ with $q^{\#}_{\ell*}(\mu)=\nu$.

The following is a slight modification of \cite[Theorem 4.1]{lia92}.
We state it here without proof.
\begin{lem}\label{perm}
For any given $\mu\in E(\f^{\#}_n, \chi^{\#})$ and any permutation
$$
r: \{ 1, 2,  \ldots , n\}\to \{r(1), r(2),  \ldots , r(n)\}
$$
there exists $\bar \mu\in E(\f^{\#}_n, \chi^{\#})$ such that  $q^{\#}_{n*}(\mu)=q^{\#}_{n*}(\bar \mu)$
and
$$
\int\,\omega_{r(i)}\,d\mu\,\,=\,\,\int\,\omega_i\,d\bar\mu,\,\,\,\,i=1,2, \ldots , n.
$$
\end{lem}

Now we fix the positive integer $\ell$, $\ell \leq n, $ as in the Main Theorem.
Define
\begin{equation}\label{idell}
id_\ell:\,\, \f^{\#}_n\to \f^{\#}_\ell : \a=(v_1,  \ldots , v_{n-\ell},
v_{n-\ell+1}, \ldots ,
v_n)\to \tilde \a=(v_{n-\ell+1},  \ldots , v_n).
\end{equation}
Then, $id_\ell$ is a continuous
projection. For $\tilde \alpha\in \f^{\#}_\ell$, set:
\begin{equation}\label{zetak}
\tilde \zeta_k(\tilde\a)
    = \zeta_{n-\ell+k} \circ (id_\ell)^{-1}(\tilde \a),
  \quad \text{ and, } \quad
\tilde \omega_k(\tilde\a)
    = \omega_{n-\ell +k} \circ (id_\ell)^{-1}(\tilde\a).
\end{equation}
It is clear by the definitions that both $\tilde \zeta_k(\tilde\a)$ and  $\tilde \omega_k
(\tilde\a)$ are independent of the choice of preimages in $id_\ell^{-1}(\a)$.
Thus  $\tilde \zeta_k,\,\,
\tilde \omega_k:\,\,
\f^{\#}_\ell
\to \RR$ are all  well defined.
For $\mu\in E(\f^{\#}_n, \chi^{\#})$ set  $\tilde \mu := id_{\ell *}(\mu)$. Then
$\tilde\mu\in E(\f^{\#}_\ell, \chi^{\#})$.
Take $\a=(u_1,  \ldots , u_n) \in Q_\mu(\f^{\#}_n, \chi^{\#})$.
Then by Lemma~\ref{logint},
$$
\lim_{t\to\pm\infty}\frac 1t\log
\zeta_{\a(n-\ell+k)}(t)
=\lim_{t\to\pm\infty}\frac 1t\int_0^t\om_{n-\ell+k}(\chi^{\#}_\tau(\a))d\tau
=\int_{\f^{{\#}}_{n}} \om_{n-\ell+k}\, d \mu,\,
$$
for $k = 1,2,\dots, \ell $.
Now write $\tilde \a := id_\ell(\a)=(u_{n-\ell+1},  \ldots , u_n)$. Then
$\tilde \omega_k(\tilde\a)=\frac{d\tilde\zeta_{\tilde\a k}(t)}{dt}\bigl|_{t=0}$ and
$\tilde\om_k(\chi^{\#}_t(\tilde\a))=\frac1{\tilde \zeta_{\tilde\a k}(t)} \frac{d\tilde \zeta_{\tilde\a k}(t)}
{dt}$
and thus Lemma~\ref{logint} holds for $\tilde \zeta_{\a k}$ and
$\tilde \omega_k$, $k=1,  \ldots , \ell$.
Observe that for $\tilde\a\in Q_{\tilde \mu}(\f^{\#}_\ell, \chi^{\#})$
we have
\begin{equation*}
\begin{split}
\lim_{t\to \infty}\frac 1t\log \tilde\zeta_{\tilde\a k}(t)
& = \int_{\f_{\ell}^{{\#}}}\tilde\omega_k d\tilde\mu\\
& = \int_{\f^{{\#}}_{n}}  \omega_{n-\ell+k}\,d \mu\\
& = \lim_{t\to \infty}\frac 1t\log \zeta_{\a(n-\ell+ k)}(t),
          \quad  k=1,  \ldots , \ell.
\end{split}
\end{equation*}

We remark that the above function $\tilde \zeta_k,\,\f^{\#}_\ell\to \RR$ is not
necessarily the same as $\zeta_k:\,\f^{\#}_\ell\to \RR$, and the function $\tilde
\omega_k:\, \f^{\#}_\ell\to \RR$ is not exactly  the same  as  $
\omega_k: \f^{\#}_\ell
\to \RR$, where  $\zeta_k,\,\,\omega_k: \f^{\#}_\ell\to \RR$
are  given in Definition 3.1 with $n$ replaced by $\ell$.

\begin{prop}\label{prop:mu}
Let $\nu \in E(M^n, \phi)$ be as in the Main Theorem,
that is, it supports $\ell$ nonzero Lyapunov exponents
$\lambda_1< \ldots <\lambda_\ell$
together with $n-\ell$ zero Lyapunov exponents.
Then there exist two probabilities
$ \mu\in E(\f^{\#}_n, \chi^{\#})$ and $\tilde \mu\in E(\f^{\#}_\ell, \chi^{\#})$,
and two subsets $\Lambda\subset M^n$ and $W\subset \f^{\#}_n$ such that
\begin{enumerate}
\item  $\,q^{\#}_{n*} (\mu)=\nu,\,\,\,\,q^{\#}_{\ell*} (\tilde\mu)=\nu,\,\,\,\, id_{\ell*}(\mu)
         =\tilde \mu;$
\item  $\phi_t(\Lambda)=\Lambda, \quad \chi_t^{\#}(W)=W, \quad
         \nu(\Lambda)=1, \quad \text{and} \quad \mu(W)=1;$
\item  For each  $ x\in \Lambda \text{ and } \a\in W \text{ with } q^{\#}_n(\a)=x $
         $$
         \lim_{t\to\pm\infty}\frac 1t\log \zeta_{\a(n-\ell+ k)}(t)
         =\int_{\f^{{\#}}_{n}}  \om_{n-\ell+k}\,  d  \mu
         =\int_{\f^{{\#}}_{\ell}} \tilde \om_k\,  d \tilde \mu
         = \lambda_k,
$$
for $k=1,2,\ldots,\ell$.
\end{enumerate}
\end{prop}

\noindent
{\bf Proof.} Take a $\phi_t$-invariant subset $\Lambda_1\subset M^n$
with $\nu$-total probability so that at each point the spectrum of all Lyapunov
exponents is  $\lambda_1,  \ldots , \lambda_\ell$ together with $n-\ell$ zeros. Furthermore,
$$
\{ \lambda_1, \ldots , \lambda_\ell\,\,\, 0, \ldots ,0 \}=
\bigl\{ \int \om_k \,  d  \mu : \mu \in
   E(\f^{\#}_n,\,\,\chi^{\#}),\,\, q^{\#}_{n*}(\mu)=\nu,\,\,k=1,2,\dots,n \bigr\}.
$$
The existence of $\Lambda_1$ follows from the hypothesis of the present
proposition and Theorem 2.2 in \cite{Sun98a}. Choose an arbitrary $\mu_1\in E(\f^{\#}_n, \chi^{\#})$
to cover $\nu$, i.e., $q^{\#}_{n*}(\mu_1)=\nu$. We claim that
\begin{equation}\label{lamom}
\{ \lambda_1,  \ldots , \lambda_\ell,\, 0, \ldots ,0\} = \bigl\{ \int \om_k\,  d  \mu_1, k =1 ,2 , \dots ,n  \bigr\}.
\end{equation}

Observe that $\mu_1(Q_{\mu_{1}}(\f^{\#}_n, \chi^{\#}))=1$ and
$\nu (q^{\#}_n Q_{\mu_{1}}(\f^{\#}_n, \chi^{\#}))=1$, thus
$$
\nu \Bigl(q^{\#}_nQ_{\mu_{1}}(\f^{\#}_n, \chi^{\#}) \bigcap \Lambda_1 \Bigr) = 1.
$$
Take $x\in q^{\#}_n
Q_{\mu_{1}}(\f^{\#}_n, \chi^{\#})\bigcap \Lambda_1$ and $\a=(u_1,  \ldots , u_n)\in \f^{\#}_n(x)
\bigcap Q_{\mu_{1}}(\f^{\#}_n, \chi^{\#})$. Remember that $\omega_k$ is
a continuous function. By Lemma~\ref{logint} we then have
$$
\lim_{t\to\pm\infty}\frac 1t\log \zeta_{\a k}(t)
=\lim_{t\to\pm\infty}\frac 1t\int_0^t\om_k(\chi^{\#}_s(\a))ds
=\int \om_k\, d \mu_1, \quad  k=1,2,\dots,n.
$$
We point out that in the  case when index $k=1$ we have
\begin{equation*}
\begin{split}
\lim_{t\to\pm\infty}\frac 1t\log \| \Phi_t(u_1) \|
& = \lim_{t\to\pm\infty}\frac 1t \log \zeta_{\a 1}(t)\\
& =\int \om_1\, d \mu_1.
\end{split}
\end{equation*}

If we suppose that Equation (\ref{lamom}) is not true,
then there would exist a minimal index $i_0>1$ such that
$$
\lim_{t\to\pm\infty}\frac 1t\log
\lV \Phi_t(u_{i_{0}})\rV \neq
\int \om_k\, d \mu_1, \quad \text{ for all }\,\, k=1,2,\dots,n.
$$
Note that $\{\proj_1\chi^{\#}_t(\a),  \ldots , \proj_n \chi_t^{\#}(\a)\}$ is an orthonormal
frame on the tangent space $T_{\phi_{t}(x)}M^n$ and
$<\Phi_t(u_{i_{0}}), \proj_j \chi^{\#}_t(\a)>=0$ for each  $j=i_{0}+1,  \ldots , n$.
We can represent $\frac {\Phi_t(u_{i_{0}})}{\lV\Phi_t(u_{i_{0}})\rV}$ as
$$
\frac {\Phi_t(u_{i_{0}})}{\lV\Phi_t(u_{i_{0}})\rV}=a_1(t) \proj_1\chi^{\#}_t(\a)+a_2(t) \proj_2\chi^{\#}_t(\a)
+ \ldots + a_{i_{0}}(t) \proj_{i_{0}}\chi^{\#}_t(\a),
$$
where $|a_k(t)| \leq 1, k=1, 2,  \ldots , i_0$.  Now let us suppose that
$\lim _{t\to \infty}|a_{i_{0}}(t)|>0$. Observe that both
$a_{i_{0}}(t) \lV \Phi_t(u_{i_{0}})\rV \proj_{i_{0}}\chi^{\#}_t(\a)$ and
$\proj_{i_{0}}\chi_t(\a)$ express the same projection of $\Phi_t(u_{i_{0}})$
on the direction determined by $\proj_{i_{0}}\chi^{\#}_t(\a)$, thus
$$
|a_{i_{0}}(t)| \lV\Phi_t(u_{i_{0}})\rV =\zeta_{\a_{ i_{0}}}(t).
$$
Therefore by Lemma~\ref{logint}
\begin{equation*}
\begin{split}
\lim_{t\to\infty}\frac 1t\log \lV \Phi_t(u_{i_{0}})\rV
& = \limsup_{t\to\infty}\frac 1t\log \lV \Phi_t(u_{i_{0}})\rV\\
& = \limsup_{t\to\infty}\frac 1t\log |a_{i_{0}}(t)|^{-1} +
\limsup_{t\to\infty}\frac 1t\log \zeta_{\a i_{0}}(t)\\
& = \lim_{t\to\infty}\frac 1t\log \zeta_{\a i_{0}}(t)\\
& = \lim_{t\to\infty}\frac 1t \int_0^t\, \omega_{i_{0}}(\chi^{\#}_s(\a))\, ds\\
& = \int \omega_{i_{0}}\, d\,\mu_1.
\end{split}
\end{equation*}
This is a contradiction to the choice of $i_0$. For the case $\lim _{t\to \infty}|a_{i_{0}}(t)|=0$,
one then gets that
$$
\lim_{t\to\infty}\frac 1t\log \lV \Phi_t(u_{i_{0}})\rV
$$
coincides with $\int\, \omega_i\,d\mu_1$ for some $i <i_0$,  again a contradiction to the
choice of $i_0$. Consequently, (\ref{lamom}) holds.

Now there is a permutation $r:\{1, 2,  \ldots , n\}\to \{ r(1), r(2),  \ldots , r(n)\}$ so that
$$
\int \omega_{r(i)}\,d\mu_1=0,\,\,\,\, i=1,2, \ldots , n-\ell,
$$
and
$$
\int \omega_{r(i)} \, d\mu_1 = \lambda_{i- (n - \ell)}, \quad i = n - \ell +1, n - \ell +2, \ldots , n.
$$
From Lemma~\ref{perm}, there exists a covering probability $\mu\in E(\f^{\#}_n,\chi^{\#})$
of $\nu$, $q^{\#}_{n*}(\mu)= q^{\#}_{n*}(\mu_1)=\nu$,
so that
$\int \omega_i\,d\mu =\int \omega_{r(i)}\,d\mu_1, \,\, \,\,i=1,2, \ldots , n$.

Define $W:=Q_\mu(\f^{\#}_n, \chi^{\#})$ and $\Lambda:= \Lambda_1\bigcap q^{\#}_n(W)$.
Then $\nu(\Lambda)=\mu(W)=1$,
$\chi^{\#}_t(W)=W,\,\, \phi_t(\Lambda)=\Lambda,\,\, t\in \RR$, and
$ q^{\#}_n(W)=\Lambda$. Define $\tilde \mu:= id_{\ell*}(\mu)$. Then
$\tilde \mu\in E(\f^{\#}_\ell, \chi^{\#})$. Clearly $q^{\#}_{\ell *}(\tilde \mu)=\nu$,
and
$$
\int \tilde\omega_k\,d\tilde \mu= \int \omega_{n-\ell +k}\, d \mu,\,\,\,\, k=1, \ldots , \ell.
$$
Take $x\in \Lambda$ and $\a\in W\bigcap \f^{\#}_n(x)$, then
$$
\lim_{t\to\pm\infty}\frac 1t\log
\zeta_{\a(n-\ell+ k)}(t)
=\int \om_{n-\ell+ k}\,  d \mu\,\,  =
\int \tilde \om_k\,  d \tilde \mu,\,\,  \,\, k=1,2,\dots,\ell.
$$
This completes the proof of Proposition~\ref{prop:mu}.

\begin{cor}
For any given $\tilde{\alpha}=(u_{n-\ell +1}, \ldots , u_n)\in
Q_{\tilde \mu}(\f^{\#}_\ell, \chi^{\#})$ with $q^{\#}_{\ell*}(\tilde \a)
\in \Lambda$ we have
$$
\lim_{t\to -\infty}\frac 1t\log \lV \Phi_t(u_i)\rV<0, \qquad
i= n-\ell+1,  \ldots , n-\ell +p,
$$
and,
$$
\lim_{t\to +\infty}\frac 1t\log \lV\Phi_t(u_i)\rV>0, \qquad i=
n-\ell + p+1, \ldots , n,
$$
where $p$ satisfies:
$\lambda_1< \ldots <\lambda_p<0<\lambda_{p+1}< \ldots <\lambda_\ell$.
\end{cor}

\noindent
{\bf Proof.}
For $n - \ell + p + 1 \le i_0 \le n $, we
have by Lemma~\ref{logint} that,
\begin{equation*}
\begin{split}
0  & < \int_{\f^{{\#}}_{n}}\omega_{i_{0}}d\mu\\
& = \lim_{t\to +\infty}\frac 1t\int_0^t\omega_{i_{0}}(\chi^{\#}_\tau(\a))dt\\
& = \lim_{t\to +\infty}\frac 1t\log\zeta_{\a i_{0}}(t)\\
& \leq \lim_{t\to \infty}\frac 1t\log\lV \Phi_t(i_{k_{0}})\rV .
\end{split}
\end{equation*}

For $n - \ell + p +1 \le i_0 \le n $, we may deduce a
similar inequality. \hfill $\Box$

%%%%%%%%%%%%%%%%%%%%%%%%%%%%%%%%%%%%%%%%%
%%%%%%%%%%%%%%%%%%%%%%%%%%%%%%%%%%%%%%%%%
%%%%%%%%%%%%%%%%%%%%%%%%%%%%%%%%%%%%%%%%%
%%%%%%%%%%%%%%%%%%%%%%%%%%%%%%%%%%%%%%%%%
%%%%%%%%%%%%%%%%%%%%%%%%%%%%%%%%%%%%%%%%%

\section{Reduced standard linear systems of $\ell$ differential equations }

We start this section from the  $\phi$-invariant, ergodic probability
$\nu \in E(M^n,\phi)$ assumed  in the Main Theorem  together with
its two covering probabilities
$\mu\in$  $ E(\f^{\#}_n$,  $ \chi^{\#})$ and  $ \tilde\mu \in $ $ E(\f^{\#}_\ell, \chi^{\#})$
and the corresponding total probability
subsets $\Lambda\subset M^n$
and $W\subset \f^{\#}_n$ as in Proposition~\ref{prop:mu}.
Take a  point $x\in \Lambda$ and an orthonormal frame $\a\in W\bigcap \f^{\#}_n(x)$.
Then

\begin{equation*}
\begin{split}
\lim_{t\to \infty}\frac 1t\log\zeta_{\a k}(t)
=&\lim_{t\to \infty}\frac 1t  \int_0^t\omega_k(\chi^{\#}_\tau(\a))d\tau\\
=&\int_{\f^{{\#}}_{n}} \omega_k d\mu\\
=&\int_{\f^{{\#}}_{\ell}}  \tilde \omega_{k-(n-\ell)}\,d\tilde \mu, \quad k= n-\ell+1,  \ldots , n.
\end{split}
\end{equation*}

In this section we will construct the reduced standard   linear system needed in the Main
Theorem  along the orbit $\mathrm{orb}(x, \phi)$ with respect to the given orthonormal frame
$\a\in \f^{\#}_n(x)$, by developing the technique in \cite{lia63}.

Since $\chi^{\#}_t(\a)$ is an orthonormal frame at $T_{\phi_t(x)}M^n$, there
exists an $n\times n$ matrix $B_\a(t)$ such that
$\Phi_t(\a)= \chi_t^{\#}(\a)\circ B_\a(t)$. Define
$R_{\a(t)}=\frac {dB_{\a}(t)}{dt} \circ B_\a(t)^{-1}$.
Define a diagonal matrix
$$
\zeta_\a(t)= \mathrm{diag} \bigl(\zeta _{\alpha 1}(t), \zeta _{\alpha 2}(t),
 \ldots ,\zeta _{\alpha n}(t) \bigr).
$$
From Gram-Schmidt orthogonalization,
$\chi_t^{\#}(\alpha)=\Phi_t(\a)\circ \Gamma(\Phi_t(\a))\circ \zeta_\a^{-1} (t)$, or,
$\Phi_t(\a)=\chi_t^{\#} \circ \zeta_\a(t)\circ \Gamma(\Phi_t(\a))^{-1}$,
where $\Gamma(\Phi_t(\a))$ is an $n\times n$ upper triangular matrix
with elements 1 on the diagonal.
So $B_\a(t)=\zeta_\a(t)\circ \Gamma(\Phi_t(\a))^{-1}$, which is
differentiable with respect to $t\in \RR$. Observe
$$
\frac 1{\zeta_{\alpha k}(t)}\frac {d\zeta _{\alpha k}(t)}{dt}=
\om_k(\chi^{\#}_t(\a)),\,\, k=1, \ldots , n ,
$$
and
$$
\frac {dB_\a(t)}{dt}\circ  B_\a(t)^{-1}=
\left (\begin{array}{cccc}
\frac 1{\zeta_{\alpha 1}(t)}\frac {d\zeta _{\alpha 1}(t)}{dt}& &* &\\
&\frac 1{\zeta _{\alpha 2}(t)}\frac {d\zeta _{\alpha 2}(t)}{dt}& &\\
&&\ddots &\\
&&&\frac 1{\zeta _{\alpha n}(t)}\frac {d\zeta _{\alpha n}(t)}{dt}
\end{array}
\right).
$$
Thus
\begin{equation}\label{Ral}
R_\alpha (t) =\left( \begin{array}{cccc}
\omega _1(\chi^{\#}_t(\a))& &* & \\
&\omega _2(\chi^{\#}_t(\a))& & \\
&&\ddots &\\
&&&\omega _{n}(\chi^{\#}_t(\a))
\end{array}
\right ).
\end{equation}

Now let us denote  $R_\a(t)^T$ by
$(r_{ij}(t))_{n\times n}$, where $r_{ij}(t)=0$ if $i<j;$ $r_{ii}(t)=
\omega_i(\chi^{\#}_t(\a))$, $i, j=1, \ldots ,n$.
Set $\tilde \a:=id_\ell(\a)$.  Recall  from Section 3  that, $\tilde \omega_i(\tilde \a)=\omega_{(n-\ell)+i}\circ
id_\ell^{-1}(\tilde\a)$, for $i=1,  \ldots , \ell$.
We define a triangular $\ell\times \ell$ matrix $A_{\ell\times \ell}(t)
=(a_{ij}(t))_{\ell\times \ell}$ as follows:
$a_{ij}(t)=0$ if $i<j;$ $a_{ij}(t)=r_{(n-\ell+i)(n-\ell+j)}(t)$ if $i>j;$
$a_{ii}(t)=\tilde \omega_{i}(\chi^{\#}_t(\tilde\a))$, $i, j =1,  \ldots , \ell$.
Thus
\begin{equation}\label{Aellell}
A_{\ell\times \ell} (t) =\left( \begin{array}{cccc}
\tilde\omega _1(\chi^{\#}_t(\tilde \a))& & & \\
&\tilde \omega _2(\chi^{\#}_t(\tilde \a))& & \\
&&\ddots &\\
&*&&\tilde \omega _{\ell}(\chi^{\#}_t(\tilde \a))
\end{array}
\right ),
\end{equation}
where $\a\in W\bigcap \f^{\#}_n(x)$ and $x\in \Lambda$.

\begin{defn}\label{def:rsls}
We call
\begin{equation}\label{RSLS}
\frac {dy}{dt}=y A_{\ell\times \ell}(t),
\end{equation}
the reduced standard linear system of  $\ell$ differential equations for  the
given system $(M^n, S, \nu)$ with respect to an orthonormal $n$-frame $\a$,
where $A_{\ell \times \ell}(t)$ is given by (\ref{Aellell}).
\end{defn}

\noindent
{\bf Proof of the Main Theorem (1.)(2.).}
For (1) it is sufficient to show $A_{\ell\times \ell}(t)$ is uniformly bounded.
In \cite{lia63} Liao proved that $\sup_{t\in \RR} \lV R_\a(t)\rV <\infty$, from which it is easy to get
$$
\sup_{t\in \RR} \lV A_{\ell\times \ell}(t)\rV \leq \sup_{t\in \RR} \lV R_\a(t)\rV <\infty.
$$
Now we prove the Main Theorem (2.) by showing the following proposition.

\begin{prop}\label{prop:vi}
Let $\nu \in E(M^n, \phi)$ be as in the Main Theorem. Let us take covering probabilities
$\mu\in E(\f^{\#}_n, \chi^{\#})$, and $\tilde \mu\in E(\f^{\#}_\ell, \chi^{\#})$, satisfying
$q^{\#}_{n*}(\mu)= \, q^{\#}_{\ell*}(\tilde\mu)=\nu$, and take a $\mu$-total
probability subset $W\subset \f^{\#}_n$ and a $\nu$-total probability subset
$\Lambda\subset M^n$ as in Proposition~\ref{prop:mu}. Take $x\in \Lambda$ and $\a
\in W\bigcap \f^{\#}_n(x)$ and construct the reduced standard linear system
(\ref{RSLS}) of $\ell$ differential equations  as in Definition~\ref{def:rsls}.
For a coordinate vector $e_i=(0,  \ldots , 0, 1(i), 0,  \ldots ,0)\in \RR^\ell$ denote by
$\tilde y(t, e_i)$ a unique solution of the initial value problem
(\ref{rsls:ivp}) with $y(0,e_i) = e_i$. Then
$$
\lim_{t\to \pm \infty}\frac 1t\log \lV \tilde y(t, e_i)\rV= \lambda_i, \quad i=1, \ldots , \ell.
$$
\end{prop}

\noindent
{\bf Proof.} Solving the  initial value problem
$$
\frac {dy_\ell}{dt}=\tilde \omega_\ell(\chi^{\#}_t(\tilde \alpha))y_\ell, \qquad
y(0)= e_\ell,
$$
we get
$$
y_\ell(t,e_\ell)= e_\ell \,e^{\int_{0}^{t}\tilde \omega_\ell(\chi^{\#}_\tau(\tilde \a))d\tau}.
$$
This equality together with  Proposition~\ref{prop:mu}  implies the following
\begin{equation*}
\begin{split}
\lim_{t\to  \infty}\frac 1t\log |  y_\ell(t, e_\ell)|
& = \lim_{t\to  \infty}\frac 1t \int _0^t\tilde \omega_\ell(\chi_\tau^{\#}(\tilde \a))d\tau\\
& = \lim_{t\to  \infty}\frac 1t \int _0^t\omega_n(\chi_\tau^{\#}(\a))d\tau\\
& = \int _{\f^{{\#}}_{n}}\omega_n \,d\mu\\
& = \int _{\f^{{\#}}_{\ell}}\tilde \omega_\ell \,d\tilde \mu\\
& = \lambda_\ell.
\end{split}
\end{equation*}

Solving the initial value problem
$$
\frac {dy_{\ell-1}}{dt}=\tilde\omega_{\ell-1}(\chi^{\#}_t(\tilde \alpha))y_{\ell-1}+ e_\ell r_{\ell(\ell-1)}(t)
e^{\int_{0}^{t}\tilde \omega_\ell(\chi^{\#}_\tau(\tilde \a))d\tau}, \quad
y_{\ell-1}(0,e_{\ell-1})= e_{\ell-1}
$$
we get
\begin{equation*}
\begin{split}
y_{\ell-1}(t, e_{\ell-1})
& = e_{\ell-1}e^{\int_{0}^{t}\tilde \omega_{\ell-1}(\chi^{\#}_\tau(\tilde \a)) d\tau}\\
& \hspace*{.3cm} +  e_\ell e^{\int_{0}^{t}\tilde \omega_{\ell-1}(\chi^{\#}_\tau(\tilde \a))d\tau}\int_0^t
r_{\ell(\ell-1)}(\tau) e^{\int_{0}^{\tau}\tilde \omega_{\ell}(\chi^{\#}_s(\tilde \a)) ds}
e^{-\int_{0}^{\tau}
\tilde \omega_{\ell-1}(\chi^{\#}_s(\tilde \a))ds}d\tau.
\end{split}
\end{equation*}
Thus,
\begin{equation*}
\begin{split}
\lim_{t\to  \infty}\frac 1t\log |  y_{\ell-1}(t, e_\ell)|
& = \lim_{t\to \infty}\frac 1t \int _0^t\tilde \omega_{\ell-1}(\chi_\tau^{\#}(\tilde \a))d\tau\\
& = \lim_{t\to \infty}\frac 1t \int _0^t\omega_{n-1}(\chi_\tau^{\#}(\a))d\tau\\
& = \int _{\f^{{\#}}_{n}}\omega_{n-1} \,d\mu\\
& = \int _{\f^{{\#}}_{\ell}}\tilde \omega_{\ell-1} \,d\tilde \mu\\
& = \lambda_{\ell-1}.
\end{split}
\end{equation*}

By repeating this procedure we will obtain:
$$
\lim_{t\to  \infty}\frac 1t\log |  y_{j}(t, e_j)| =\lambda_{j},
$$
for $j=1,  \ldots , \ell$.

From the form of the functions $y_\ell(t, e_\ell), \ldots, y_1(t, e_1)$,
which depend linearly on the initial values $e_1$,  $e_2$, $ \ldots $,  $e_\ell$,
we get easily
$$
\lV \tilde y(t, e_i) \rV= |y_i(t, e_i)|.
$$
Therefore
\begin{equation*}
\begin{split}
\lim_{t\to \infty}\frac{1}{t} \log \lV \tilde y(t, e_i)\rV
& = \lim_{t\to \infty} \frac{1}{t} \log | y_i(t, e_i)|\\
& = \lambda_i.
\end{split}
\end{equation*}
This  proves  the proposition and thus  proves parts (1.) and (2.) of the Main Theorem.

%%%%%%%%%%%%%%%%%%%%%%%%%%%%%%%%%%%
%%%%%%%%%%%%%%%%%%%%%%%%%%%%%%%%%%%
%%%%%%%%%%%%%%%%%%%%%%%%%%%%%%%%%%%
%%%%%%%%%%%%%%%%%%%%%%%%%%%%%%%%%%%
%%%%%%%%%%%%%%%%%%%%%%%%%%%%%%%%%%%

\section{Proof of the Main Theorem (3.)}

In this section we will complete the proof of  the Main Theorem.

Let  $\nu\in E(M^n,\phi)$ denote the  given probability in  the Main Theorem.
Let $\mu\in E(\f^{\#}_n,\chi^{\#})$
and $\tilde \mu\in E(\f^{\#}_\ell, \chi^{\#})$ be the covering probabilities,
and let $\Lambda\subset M^n$
and $W\subset \f^{\#}_n$ be the two total probability sets, $q^{\#}_n(W)=\Lambda$, as in
Proposition~\ref{prop:mu}.
Write
$\vartheta_i(\tilde \mu)=\int_{\f^{{\#}}_\ell }\,\tilde\omega_i\,d\tilde
\mu$, $i=1, 2, \ldots ,\ell$.
Then we have that $\vartheta$ satisfies:
$$
 \vartheta_1(\tilde \mu)=\lambda_1<\vartheta_2(\tilde \mu)=\lambda_2< \ldots
<\vartheta_\ell(\tilde \mu)=\lambda_\ell.
$$
Let $T_1\geq 1 $ be a fixed constant and let $ T_{i+1}=2T_i,\,\, i=1, 2,\dots$.

Recall from Section 3 the projection map (\ref{idell}).
\begin{defn}
For $\eta>0$ we denote
by $ D(\vartheta, \eta)$ the set of all $\tilde \gamma\in id_\ell(W)$ with the
property that for each integer $i\geq 1$ there exist an integer $c=c(\tilde \gamma, i,
 \eta)
\geq i$ and a sequence
\begin{equation}\label{seq}
\begin{split}
\ldots & < s(-2) < s(-1) <s(0) = 0 < s(1) < s(2) < \ldots \\
        & \lim_{j\to -\infty}s(j) = -\infty,   \qquad   \lim_{j\to +\infty} s(j) = +\infty,
\end{split}
\end{equation}
such that
$$
\frac 1l\sum_{\tau=0}^{l-1}\max_{k=1,2, \ldots , \ell}
  \Bigl| \vartheta_k(\tilde\mu)-
\frac 1{\delta T_{c}}\int_{\tau\delta T_{c}}^{(\tau+1)\delta T_{c}}\,
\tilde\omega_k(\chi^{\#}_{t+s(j)T_{c}}(\tilde\gamma))dt \Bigr| < \eta,
$$
$$
l=1, 2,  \ldots ;\,\, j=0,\pm 1, \pm 2,  \ldots ; \,\,\delta=\pm 1.
$$
\end{defn}

\begin{lem}\label{lem:mut}
$\tilde\mu( D(\vartheta, \eta))>0$.
\end{lem}
\noindent
{\bf Proof.} This is a partial  result of \cite[Theorem 2.1]{lia97}, where Liao gave a general proof.
We present a proof of our case here for convenience to readers. Set
\begin{equation*}
\begin{split}
h_k(\tilde\gamma, T, \delta)&=|\vartheta_k(\tilde\mu)-
\frac 1{\delta T}\int_0^{\delta T}\,
\tilde\omega_k(\chi^{\#}_{t}(\tilde\gamma))dt|, \quad
        0<T<+\infty,\,\,\delta=\pm 1,\,\,\tilde\gamma\in \f^{\#}_\ell,  \\
h(\tilde\gamma, T, \delta)&=\max_{k=1,2, \ldots ,\ell}h_k(\tilde\gamma, T, \delta).
\end{split}
\end{equation*}
When  $\tilde\gamma\in id_\ell(W)$ we  choose and fix $\gamma\in W$ with
$id_\ell(\gamma) =\tilde \gamma$.
We have from Proposition~\ref{prop:mu}, for $k=1, 2, \ldots ,\ell$,
\begin{equation*}
\begin{split}
\lim_{T\to \pm \infty}\frac 1T\int_0^T \tilde\omega_k(\chi_t^{\#}(\tilde \gamma))dt
& = \lim_{T\to  \pm\infty}\frac 1T\int_0^T \omega_{n-\ell+k}(\chi_t^{\#}( \gamma))dt\\
& = \vartheta_k(\tilde \mu)\\
& = \lambda_k.
\end{split}
\end{equation*}
By Chapter 6 in \cite{NS},
$$
\lim_{T\to \infty}\int_{\f^{\#}_{\ell}} \Bigl|\vartheta_k(\tilde \mu)-
\frac 1{\delta T}\int_{a}^{a+T}\,
\tilde\omega_k(\chi^{\#}_{t}(\tilde\gamma))dt \Bigr| d\tilde \mu \,=\, 0,
$$
where the convergence is uniform with respect to the choice of $a\in \RR$.
For $\delta=1$  and $a=0$ we get
$$
\lim_{T\to \infty}\int_{\f^{\#}_{\ell}}\,h_k(\tilde \gamma, T, +1)\,d\tilde \mu
=\, 0,\,\,\,\, k=1,  \ldots , \ell.
$$
For $\delta =-1$, taking $a=-T$, we get
$$
\lim_{T\to \infty}\int_{\f^{\#}_{\ell}} h_k(\tilde \gamma, T, -1)\,d\tilde \mu =\, 0, \quad
               k=1, \ldots , \ell.
$$
Therefore
$$
\lim_{T\to \infty}\int_{\f^{\#}_{\ell}} h(\tilde\gamma, T, \delta)d\tilde \mu=0,
\quad \delta =\pm 1.
$$
For $\eta>0$ one can thus take an integer $d=d(\eta)>0$ such that
$$
\int_{\f^{\#}_{\ell}} h(\tilde \gamma, T_d, \delta)d\tilde\mu
    < \frac \eta{30},\,\,\,\,\delta=\pm 1.
$$

Let us consider a $\tilde \mu$-preserving homeomorphism
$\rho=\chi_{T_{d}}^{\#}: \,\,\f_\ell^{\#}\to \f^{\#}_\ell$.
Applying the Birkhoff Ergodic Theorem to the homeomorphism $\rho^{\delta}$ and
continuous function $h(\tilde\gamma; T, \delta), \,\,\delta=\pm 1$, there is a
$\tilde \mu$-measurable
subset $X\subset \f^{\#}_\ell$ with $\tilde\mu(X)=1$ such that for any
$\tilde\gamma\in X,\,\, \delta=\pm 1$,
the following limit exists
$$
\lim_{l\to \infty}\frac 1l \sum_{\tau=0}^{l-1}h( \rho^{\delta \tau}(\tilde\gamma),
T_d,\delta)=h^*(\tilde\gamma,T_d, \delta).
$$
Moreover,
$$
\int_{\f^{\#}_{\ell}} h^*(\tilde\gamma, T_d, \delta)d\tilde \mu
=\int_{\f^{\#}_{\ell}} h(\tilde\gamma, T_d, \delta)d \tilde \mu < \frac \eta{30}.
$$
This implies that the set
$\{\tilde\gamma\in \f_\ell^{\#}\,|\, h^*(\tilde\gamma, T_d,\delta)>\frac \eta 2\}$
is $\tilde \mu$-measurable and has $\tilde\mu$-probability less than or equal to
$\frac 1{12},\,\,\delta=\pm 1$. Applying Egoroff's Theorem (see, for example, [3]) there exists
a subset $Y$ of $X$ with $\tilde \mu(Y)\geq \frac 34>0$ such that
$$
\frac 1l \sum_{\tau=0}^{l-1} h(\rho^{\delta \tau}(\tilde\gamma), T_d,\delta) <\eta, \quad \
                      \forall \tilde\gamma\in Y,\,\, l\geq \bar l.
$$
Therefore
$$
\frac 1l \sum_{\tau=0}^{l-1}\max_{k=1,2, \ldots ,\ell}
      \Bigl| \vartheta_k(\tilde \mu)-
           \frac 1{\delta T_{d}}\int_0^{\delta T_{d}}\tilde\omega_k
      (\chi^{\#}_{t+\delta \tau T_{d}}(\tilde\gamma))dt \Bigr| < \eta,
$$
or
$$
\frac 1l \sum_{\tau=0}^{l-1}\max_{k=1,2, \ldots ,\ell}\bigl|\vartheta_k(\tilde
\mu)- \frac 1{\delta T_{d}}\int_{\delta \tau T_{d}}^{\delta (\tau+1)T_{d}}
\tilde\omega_k (\chi^{\#}_t
(\tilde\gamma))dt \bigr| < \eta,\,\,\forall \tilde\gamma\in Y,\,\,\,\,l\geq \bar l.
$$
From the Poincar\'{e} Recurrence Theorem, let us take  a subset $Y'$ of $Y$,
$\tilde \mu(Y')=\tilde \mu(Y)$,  with the property that
for each  $\tilde\gamma\in Y'$ there exists a sequence $\{s(j)\}$ of the form (\ref{seq}),
so that $\rho^{s(j)}(\tilde\gamma)\in Y$. This gives rise to
$$
\frac 1l \sum_{\tau=0}^{l-1}\max_{k=1,2, \ldots ,\ell} \bigl|\vartheta_k(\tilde\mu)-
\frac 1{\delta T_{d}}\int_{\delta \tau T_{d}}^{\delta (\tau+1)T_{d}}
\tilde\omega_k
(\chi_{t+s(j)T_{d}}(\tilde\gamma))dt \bigr|< \eta,
$$
$$
l\geq\bar l, \quad  j=0, \pm 1,\pm 2, \ldots, \quad\delta=\pm 1.
$$

Denote by $\xi(\gamma)$ a character function for $Y$ on $\f_\ell^{\#}$.
Let us consider a $\tilde\mu$-preserving homeomorphism
$$
\psi: \f_\ell^{\#}\to\f^{\#}_\ell.
$$
Set
$$
\bar\xi(\tilde\gamma, \psi, \delta):=\limsup_{l\to \infty}
\frac 1l \sum_{\tau=0}^{l-1}\xi(\psi^{\delta \tau}(\tilde\gamma)), \quad
\tilde\gamma\in
\f_\ell^{\#},\,\,\delta=\pm1.
$$
Then $\bar \xi$ is a Baire function. Let
$$
E(\eta, \psi, \delta)=\{\tilde\gamma\in \f^{\#}_\ell\,|\, \bar\xi(\tilde\gamma, \psi,
\delta)>0\},
$$
$$
E(\eta, \psi)=E(\eta, \psi, -1)\bigcap E(\eta, \psi, +1).
$$
By the  Birkhoff Ergodic Theorem  there exists a subset $Z\subset \f^{\#}_\ell$,
$\tilde \mu(Z)=1$, such that for all $\tilde\gamma\in Y\bigcap Z$,
the limit exists
$$
\lim_{l\to \infty}\frac 1l \sum_{\tau=0}^{l-1}\xi(\psi^{\delta \tau}
(\tilde\gamma))=\xi^*(\tilde\gamma, \psi,\delta).
$$
Since
$$
1\geq \bar \xi(\tilde \gamma, \psi, \delta)\geq \xi^*(\tilde \gamma, \psi, \delta)\geq 0,
$$
then
\begin{equation*}
\begin{split}
\mu(E(\eta, \psi, \delta))
&\geq \int \bar\xi(\tilde\gamma, \psi, \delta)d\tilde \mu\\
&\geq \int \xi^*(\tilde\gamma, \psi, \delta)d\tilde \mu\\
&=\tilde\mu (Y\bigcap Z)\\
&=\tilde\mu(Y)\\
&\geq \frac 34,
\end{split}
\end{equation*}
for both $\delta =1$ and $\delta =-1$, which implies then
$\tilde\mu(E(\eta, \psi))\geq \frac 12$.

Now for each integer $i\geq 1$ take $\psi$ as
$$
\psi_i=\chi^{\#}_{T_{c}}
$$
where $T=T_c,\,\, c=c(i, \eta)\geq 1$ with $2^{c(i, \eta)-d}\geq \bar l$.
Moreover we take $c(i, \eta)<c(i+1, \eta),\,\, i=1,2,3,\dots$.
Then $T_{c(i, \eta)}=2^{c(i, \eta)-d}T_d\geq \bar l T_d$.
Write $F(\eta, \psi_i)=Y'\bigcap E(\eta, \psi_i)$ then
$\tilde \mu(F(\eta, \psi_i)) \geq \frac 14$.
Take $\tilde\gamma\in F(\eta, \psi_i)$. Then there is a
sequence $\{s(j)\}$ of the form (\ref{seq})
such that $\xi(\psi_i^{s(j)}(\tilde\gamma))=1$, namely, $\psi_i^{s(j)}(\tilde\gamma)
\in Y,\,\, j=0,\pm 1,\pm 2, \ldots $.  Recall by definition $T_c=2^{c-d} T_d$. We have
\begin{equation*}
\begin{split}
\frac 1l \sum_{\tau=0}^{l-1}\max_{k=1,2, \ldots ,\ell} & \Bigl|\vartheta_k(\mu)-
\frac 1{\delta T_{c}}\int_{\delta \tau T_{c}}^{\delta (\tau+1)T_{c}}\tilde
\omega_k
(\chi^{\#}_{t+s(j)T_{c}}(\tilde\gamma))dt \Bigr| \\
& = \frac 1l \sum_{\tau=0}^{l-1}\max_{k=1,2, \ldots ,\ell} \Bigl| \vartheta_k(\mu)-
\frac 1{\delta T_{c}}\int_{\delta \tau T_{c}}^{\delta (\tau+1)T_{c}}\tilde
\omega_k
(\chi^{\#}_{t}(\psi_i^{s(j)}( \tilde\gamma)))dt \Bigr| \\
&\leq
\frac 1{l 2^{c-d}} \sum_{\tau=0}^{l 2^{c-d}-1}\max_{k=1,2, \ldots ,\ell}
\Bigl| \vartheta_k(\mu)-
\frac 1{\delta T_{d}}\int_{\delta \tau T_{d}}^{\delta (\tau+1)T_{d}}
\tilde\omega_k
(\chi^{\#}_{t}(\psi_i^{s(j)}(\tilde\gamma)))dt \Bigr|  \\
& = \frac 1{l 2^{c-d}} \sum_{\tau=0}^{l 2^{c-d}-1}\max_{k=1,2, \ldots ,\ell}
\Bigl|\vartheta_k(\mu)-
\frac 1{\delta T_{d}}\int_{0 }^{\delta T_{d}}
\tilde\omega_k
(\chi^{\#}_{t} (\rho^{\delta \tau }(\psi_i^{s(j)}(\tilde\gamma)))dt \Bigr| \\
& = \frac 1{l 2^{c-d}} \sum_{\tau=0}^{l 2^{c-d}-1}
h(\rho^{\delta \tau}(\psi_i^{s(j)}(\tilde\gamma)), T_d, \delta)\\
& < \eta,
\end{split}
\end{equation*}
where
$\tilde \gamma\in F(\eta, \psi_i);\,\,l=1, 2, \ldots$, $j=0,\pm 1,
\pm 2, \ldots$,  and  $\delta=\pm 1$.
Letting $ D(\vartheta, \eta):= F(\eta, \psi_i)$ we complete the proof of Lemma 5.1.
\hfill $\Box$

\begin{cor}\label{cor:mut}
Set $F(\eta)=\bigcap_{i=1, 2,  \ldots }F(\eta, \psi_i)$,
where $F(\eta, \psi_i)$ is as in the proof of Lemma~\ref{lem:mut}. Then $F(\eta)$ is
a Borel subset. Since $\tilde \mu(F(\eta, \psi_i))\geq \frac 14$
and  $c(i, \eta)<c(i+1, \eta)$, and thus
$F(\eta, \psi_i)\supset F(\eta, \psi_{i+1}),\,\,i=1,2, \ldots $, we then have
$$
\mu(F(\eta))\geq \frac{1}{4}>0.
$$
\end{cor}

\begin{thm}[Liao, \cite{lia97}]\label{liao8}
Consider two systems
\begin{equation}\label{liaode}
\frac {dy}{dt}=y C(t)+f(t, y), \quad (t, y) \in \RR \times \RR^\ell
\end{equation}
\begin{equation}\label{liaoiv}
\frac {dy}{dt}=y C(t)\,\,\,\,\,\,(t, y)\in \RR\times \RR^\ell.
\end{equation}
Let the following (i)(ii) and (iii) hold.

(i). For any $t\in \RR$, $C(t)=(c_{ij})_{\ell\times \ell}$ is a lower triangular $\ell\times\ell$ matrix.
$C(t)$ is continuous with respect to $t$ and  uniformly bounded.

(ii). There exist constants $\lambda>0$, $T>0$, $c=\frac 1{16}\min\{ 1, \lambda\}$,
and a bi-infinite sequence $\{s(j)\}$ of the form (\ref{seq})
so that for some integer $p\in <0, \ell>$ the following inequalities hold
$$
-\lambda <\frac 1l \sum_{\tau=0}^{l-1}\max\{-\lambda,
\max_{k=1,2, \ldots ,p}
\frac 1{ \delta T}\int_{\delta \tau  T}^{\delta (\tau+1) T}
c_{kk}(t+s(j)T)dt\}< -\lambda+c,
$$
$$
\lambda-c<\frac 1l \sum_{\tau=0}^{l-1}\min\{\lambda,
\min_{k=p+1, \ldots , \ell}
\frac 1{\delta T}\int_{\delta \tau T}^{\delta (\tau+1)T}
c_{kk}(t+s(j)T)dt\}< \lambda
$$
$$
j=0,\pm 1, \pm 2,  \ldots ;\,\,l=1,2,  \ldots ;\,\,\delta=\pm1.
$$

(iii). Vector function f(t, y) is continuous with $(t, y)$ and is uniformly bounded
and  Lipschitz with respect to $y$.

Then, for each $u^*\in \RR^\ell$ there exists uniquely $u\in \RR^\ell$  so that
the solutions $y(t, u^*)$ and  $y(t, u)$ of the initial value problem (\ref{liaode}), (\ref{liaoiv})
with initial conditions $y(0;  u^*)=u^*$ and $y(0;  u)=u$,
respectively, satisfy the following relation.

(a). There is a integer sequence
$$
\ldots < m(-2) < m(-1) <m(0) = 0 < m(1) < m(2) < \ldots
$$
$$
\lim_{j\to -\infty}m(j) = -\infty,   \qquad   \lim_{j\to +\infty}m(j) = +\infty
$$
so that
$$
\sup_{k\in Z}\lV y(m(k)T, u)-y(m(k)T, u^*)\rV <\infty.
$$

(b). The map $\Delta^*: \RR^\ell\to \RR^\ell,\,\, u^* \to u$ is surjective.

(c). There exist constants $C^*>0$ and $d>0$ so that
$$
\lV y(t, \Delta^*(u^*))-y(t, u^*)\rV\leq C^* \exp(2c|t-s(j)T|+d),
    \quad  j=0, \pm 1, \pm 2, \, \ldots .
$$
\end{thm}
\noindent
{\bf Proof.} (a) and (b) are  Theorem 3.1 in \cite{lia97}, its   Corollary 1 is (c).

\noindent
{\bf Proof of  the Main Theorem (3.).}
For $\nu\in E(M^n, \phi)$ let us consider  all its $\ell$
nonzero Lyapunov exponents
$\lambda_1< \ldots <\lambda_p<\lambda_{p+1}< \ldots <\lambda_\ell$,
where $\lambda_p<0<\lambda_{p+1}$.
We recall again from Proposition~\ref{prop:mu} the covering probabilities
$\mu\in E(\f^{\#}_n,\chi^{\#}), \tilde \mu\in E(\f^{\#}_\ell, \chi^{\#})$,
$q^{\#}_{n*}(\mu)=\nu = q^{\#}_{\ell *}(\tilde \mu)$, and the subsets
$W\subset \f^{\#}_n$ and $\Lambda\subset M^n$ with
$q^{\#}_n(W)=\Lambda$. And consider continuous functions
$\tilde \zeta_{\a k},\,\,\tilde \omega_k:\,\,\f^{\#}_\ell \to \RR$ as in Section 3.
Take an arbitrary positive real $\lambda$ with  $\lambda_p< \lambda < \lambda_{p+1}$
and
$$
\lambda<\frac 12 \min_{1\leq i\neq j\leq \ell} \{|\lambda_i-\lambda_j|,\,\, |\lambda_i-0|\}.
$$
Write $c:=\frac 1{16}\min\{ 1, \, \lambda\}$ as in Theorem~\ref{liao8} and  write
$\eta:=\frac c 2$ as in Lemma~\ref{lem:mut}. We take and fix an
orthonormal $\ell$-frame
$$
\tilde \a\in F(\eta),
$$
where $F(\eta)$ is defined in the Corollary~\ref{cor:mut}. Recall by construction
$F(\eta)\subset id_\ell(W)$, one can take $\a\in W$ with $id_\ell( \a)=\tilde\a$.
By using the moving orthonormal $n$-frame
$$
\{\chi^{\#}_t(\a);\,\, t\in \RR\}
$$
we can construct as in Section 4 a reduced standard linear system (\ref{RSLS})
of $\ell$ differential equations.
As in Section 4 we can prove the Main Theorem(i)(ii) with respect to this linear system
of differential equations.

Now let us consider a  perturbed system (\ref{pert})
where $f(t, y)$ is Lipschitz and uniformly bounded.
Observe that the $kk$-th entry  of the matrix $A_{\ell\times \ell}(t)$ is
$$
a_{kk}(t)\,=\,\,\tilde \omega_k(\chi^{\#}_t(\tilde \a)), \qquad k=1, 2, \ldots , \ell.
$$
Since $\tilde \alpha \in F(\eta)$ and $p\leq \ell$ there exist, by Lemma~\ref{lem:mut} and Corollary~\ref{cor:mut}
a positive number $T>0$ and a sequence $\{s(j)\}$ of the form (\ref{seq}) such that
$$
\frac 1l\sum_{\tau=0}^{l-1}\max_{k=1,2, \ldots , p} \bigl|\vartheta_k(\tilde\mu)-
\frac 1{\delta T}\int_{\tau\delta T}^{(\tau+1)\delta T}\,
\tilde\omega_k(\chi^{\#}_{t+s(j)T}(\tilde\gamma))dt \bigr|< \eta,
$$
$$
l=1, 2,  \ldots ;\,\, j=0,\pm 1, \pm 2,  \ldots ; \,\,\delta=\pm 1.
$$
Observe $\lambda_k=\vartheta_k(\tilde \mu)<\lambda $, and so we get

\begin{equation*}
\begin{split}
\frac 1l\sum_{\tau=0}^{l-1} \Bigl( \max_{k=1,2, \ldots , p}\frac{1}{\delta T}&\int_{\tau\delta T}^{(\tau+1)\delta T}\,
\tilde\omega_k(\chi^{\#}_{t+s(j)T}(\tilde\gamma))dt\,-\,\lambda \Bigr) \\
&= \frac 1l\sum_{\tau=0}^{l-1}\max_{k=1,2, \ldots , p} \Bigl( \frac 1{\delta T}\int_{\tau\delta T}^{(\tau+1)\delta T}\,
\tilde\omega_k(\chi^{\#}_{t+s(j)T}(\tilde\gamma))dt\,-\,\lambda \Bigr) \\
&\leq \frac 1l\sum_{\tau=0}^{l-1}\max_{k=1,2, \ldots , p} \Bigr( \frac 1{\delta T}
\int_{\tau\delta T}^{(\tau+1)\delta T}\,
\tilde\omega_k(\chi^{\#}_{t+s(j)T}(\tilde\gamma))dt\,-\,\lambda_k \Bigr) \\
& \leq \frac 1l\sum_{\tau=0}^{l-1}\max_{k=1,2, \ldots , p}
\Bigr| \vartheta_k(\tilde\mu)-\frac 1{\delta T}\int_{\tau\delta T}^{(\tau+1)\delta T}\,
\tilde\omega_k(\chi^{\#}_{t+s(j)T}(\tilde\gamma))dt \Bigr|\\
&\leq \eta\\
& = \frac c 2.
\end{split}
\end{equation*}
Therefore
$$
-\lambda < \frac 1l \sum_{\tau=0}^{l-1}\max\{-\lambda, \max_{k=1,2, \ldots ,p}
\frac 1{\delta T}\int_{\delta \tau T}^{\delta (\tau+1)T}\tilde\omega_k
(\chi^{\#}_{t+s(j)T}(\gamma))dt\}< -\lambda+c,
$$
for $j=0,\pm 1, \pm 2,  \ldots ;\,\,l=1,2,  \ldots ;\,\,\delta=\pm1$.
Similarly,
$$
\lambda-c<\frac 1l \sum_{\tau=0}^{l-1}\min\{\lambda,
\min_{k=p+1,  \ldots , \ell}
\frac 1{\delta T}\int_{\delta \tau T}^{\delta (\tau+1)T}\tilde  \omega_k
(\chi^{\#}_{t+s(j)T}(\gamma))dt\}< \lambda,
$$
for $ j=0,\pm 1, \pm 2,  \ldots ;\,\,l=1,2,  \ldots ;\,\,\delta=\pm1$.

Now we apply Theorem~\ref{liao8} to complete the Main Theorem.
Since $\Delta^*$ in Theorem~\ref{liao8} is surjective,
for
$$
u_k=(0,  \ldots , 0, 1(k), 0,  \ldots , 0)\in \RR^\ell
$$
there exist
$u^*_k \in \RR^\ell$
so that $\Delta^*(u_k^*)=u_k,\,\,k= 1,  \ldots , \ell$.
From  Theorem~\ref{liao8}  the solution $y(t, u_k)$
of the  initial value problem
$$
 \frac {dy}{dt}=y A_{\ell\times\ell}(t), \qquad
          y(0, u_k)=u_k
$$
and the solution $y(t, u_k^*)$ of the  initial value problem
$$
 \frac{dy}{dt}=y A_{\ell\times \ell}(t) +f(t, y), \qquad   y(0, u_k^*)=u_k^*
$$
satisfy the relation
$$
\lV y(t, u_k^*)-y(t, u_k)\rV\leq C^*\exp(2c(|t-s(j)T|+d))
$$
for some constants $C^*>0$ and $d>0$.
Letting $j=0$ and thus $s(j)=0$ it follows
\begin{equation*}
\begin{split}
\lV y(t, u^*_k)\rV
& \leq \lV y(t, u_k)\rV \,+\,C^*\exp(2c|t|+d)\\
&\leq \lV y(t, u_k)\rV \,\times \,C^*\exp(2c|t|+d)
\end{split}
\end{equation*}
for $|t|\geq \bar t>0$. This yields by Proposition~\ref{prop:vi}
\begin{equation*}
\begin{split}
\limsup_{t\to\infty}\frac 1t\log \lV y(t, u_k^*)\rV
\leq&\limsup_{t\to \infty}\frac 1t\log \lV y(t, u_k)\rV + 2c\\
&= \lambda_k+ 2c,
\end{split}
\end{equation*}
where we recall $c=\frac 1 {16}\min\{1, \lambda\}$.
Since $\lambda$ and thus $c$ can be taken small enough, we get
$$
\limsup\frac 1t\log\lV y(t, u_k^*)\rV \leq \lambda_k, \qquad
k=1, 2,  \ldots , \ell.
$$
Now one can easily get
$$
 \lV y(t, u_k)  \rV
\leq \lV  y(t, u^*_k)\rV
\,\times \,C^*\exp(2c|t|+d)
$$
for $|t|\geq \bar t>0$.
This gives rise to
\begin{equation*}
\begin{split}
\lambda_k =&\liminf_{t\to\infty}\frac 1t\log \lV y(t, u_k)\rV\\
\leq&\liminf _{t\to \infty}\frac 1t\log \lV y(t, u_k^*)\rV + 2c.
\end{split}
\end{equation*}
Thus
\begin{equation*}
\begin{split}
\lambda_k -2c
& < \liminf _{t\to\infty}\frac 1t\log \lV y(t, u_k^*)\rV\\
& \leq \limsup_{t\to \infty}\frac 1t\log \lV y(t, u_k)\rV \\
& < \lambda_k+ 2c.
\end{split}
\end{equation*}
Since $c$ can be taken small enough, we get
$$
\lim_{t\to \infty}\frac 1t\log\lV y(t, u_k^*)\rV =\lambda_k, \,\,k=1, 2, \ldots , \ell.
$$
This completes the proof of the Main Theorem.
\hfill $\Box$

%End proof of Main Thm

\noindent
{\bf Example.} When $\ell<n-1$, the system  $(M^n,\phi, \nu)$ is not hyperbolic.
In this case   the Main Theorem does not hold for the linear system \cite[Chapter 2]{lia96}
$$
\frac {dy}{dt}=y R_{\gamma}(t)^T
$$
of $n$ first order differential equations based on $\a\in \f^{\#}_n$, where
$R_\a(t)$ is defined as in Section 4 (see also \cite[Chapter 2]{lia96}).
This is illustrated by the following
example. Let $n=2,\,\, \ell=1$. Take $ \a=( u_1, u_2)$
as in Section 3. Then the linear system based on $ \a$ is
$$
(\frac{dy_1}{dt},\frac{dy_2}{dt})=(y_1, y_2)
\left(\begin{array}{cc}
\omega _1(\chi^{\#}_t(\a)) & \\
&\omega _2(\chi^{\#}_t(\a))
\end{array}
\right).
$$
We consider the case when
$\lim_{t\to \pm\infty}\frac 1t\int_0^t\omega_t( \chi^{\#}_t(\a))dt =\lambda<0$
and $\omega_2(\chi^{\#}_t(\a))\equiv 0\,\,\forall t\in \RR$.
Let us consider a perturbed system
$$
(\frac{dy_1}{dt},\frac{dy_2}{dt})=(y_1, y_2)
\left(\begin{array}{cc}
\omega _1(\chi^{\#}_t(\a)) & \\
&\omega _2(\chi^{\#}_t(\a))
\end{array}
\right)
+\left(\begin{array}{c} a\\
                         a
\end{array}
\right ),
$$
where $a>0$ is a small constant. We get $y_2(t)=at$,  and thus get
$$
0 > \lim_{t\to \pm \infty}\frac 1t\log \lV (y_1(t), y_2(t)\rV
      \geq \lim_{t\to+\infty}\frac 1t\log|at|=0,
$$
which is a contradiction.

%%%%%%%%%%%%%%%%%%%%%%%%%%%%%%%%%%%%%%%%%
%%%%%%%%%%%%%%%%%%%%%%%%%%%%%%%%%%%%%%%%%
%%%%%%%%%%%%%%%%%%%%%%%%%%%%%%%%%%%%%%%%%

\section{A persistence property  for Liao perturbations}

A nearby $C^1$ vector field, while perturbing a given one, keeps neither
 value nor  sign of Lyapunov exponents, in general. However, if we perturb the
given $C^1$ vector field by a ``Liao perturbation", we will show in this section that
the perturbed vector field  will keep both  sign and
value of the nonzero Lyapunov exponents.
The class of Liao perturbations is constructed using the standard
system of the given vector field.

Recall that $S$ is the  $C^1$ vector field on  $M^n$ given in Section 1. It reduces
then in Section 2 the flows $\phi: \,M^n\to M^n, \quad \chi^{\#}:\, \f^{\#}_n\to \f^{\#}_n$.
Let  $\nu\in E(M^n,\phi)$ denote the  probability in the Main Theorem.
Let $\eta>0$ be small and let $F(\eta)$ be as in the Corollary~\ref{cor:mut},
$\tilde \mu(F(\eta))>\frac 14$. From Lemma~\ref{lem:mut}, $F(\eta)\subset id_\ell(W)$.
Recall from Section 3 the projection map $id_\ell: \f^{\#}_n\to \f^{\#}_\ell$.

Now we recall briefly the Liao standard system for a perturbation vector
field \cite[Chapter 2]{lia96} with respect to the  orthonormal $n$-frame $\beta$ we
chose. Let us take and fix $x\in \Lambda$ and $\beta\in W$ so that
$q_{n}^{\#}(\beta)=x$ and $\tilde \beta:= id_\ell(\beta)\in F(\eta)$.
Construct a standard map ${\mathcal P}_\beta:\, \RR\times \RR^{n}\to M^n$
$$
{\mathcal P}_\beta(t, y)= \exp (\sum_{i=1}^n \,y^i \proj_i\chi^{\#}_i(\beta)), \quad
y=(y^1,  \ldots , y^n).
$$
As $M^n$ is a compact $C^\infty$ Riemannian manifold, the exponential map
$\exp: TM^n\to M^n$ is $C^\infty$
and there exists  a constant $\zeta_0>0$ such that for any $x\in M^n$, $\exp$ maps
$\{u\in T_xM^n\,|\,\,\lV u\rV<\zeta_0\}$ differentially into a neighborhood of $x$ on $M^n$.
Let $X$ be a $C^1$ vector field, a perturbation to the given vector field
$S$.
Fixing $t\in \RR$, there exists a unique tangent vector field $X_\beta(t, y)$ on
$$
B_0=\{ y\in \RR^n\,|\,\lV y\rV <\zeta_0\}
$$
so that
\begin{equation*}
\begin{split}
d{\mathcal P}_{\beta t}(X_\beta(t, y))
& = d{\mathcal P}_{\beta }(0, X_\beta(t, y))\\
& = X({\mathcal P}_\beta(t, y))-
d{\mathcal P}_{\beta }(\frac \partial {\partial t}|_{(t, y)}).
\end{split}
\end{equation*}
The system
$$
\frac {dy}{dt}=X_\beta(t, y)
$$
can be written as
\begin{equation}\label{liaoss}
\frac {dy}{dt}=y R_\beta(t)^T +\bar f(t, y),
\end{equation}
where $R_\beta(t)^T=(r_{ij})_{n\times n}$ is defined in Section 4,
$r_{ii}(t)=\omega_i(\chi^{\#}_t(\beta))$, for $i= 1,  \ldots , n$.
The vector function $\bar f(t, y)$   is bounded and  Lipschitz.
The  system (\ref{liaoss}), called the Liao standard system of $X$ based on $(S, \nu)$,
was employed by Liao  to prove the $C^1$ closing lemma
\cite[Appendix A]{lia96} and topological stability for Anosov flows [7, Chapter 2].

Based on the Liao standard system, we now introduce the terminology of
Liao perturbation to the given vector field $(M^n, S, \nu)$ in our Main Theorem.
We define a triangular $\ell\times \ell$ matrix
$A_{\ell\times \ell}(t) =(a_{ij}(t))_{\ell\times \ell}$ as follows:
$a_{ij}(t)=0$ if $i<j;$ $a_{ij}(t)=r_{(n-\ell+i)(n-\ell+j)}(t)$ if $i>j;$
$a_{ii}(t)=\tilde \omega_{i}(\chi^{\#}_t(\tilde\beta))$, $i, j =1,  \ldots , \ell$.
And define a vector function $f: \RR^\ell\to \RR$ by
$$
f_i(t, y)=\bar f_{n-\ell+i}(t, (0,  \ldots ,0(n-\ell), y^1,  \ldots , y^\ell), \quad i=1,  \ldots , \ell.
$$
We then call the system
\begin{equation}\label{RSS}
\frac {dy}{dt} = y A_{\ell\times \ell}(t)  + f(t, y),
\end{equation}
a reduced standard system of the perturbation vector field $X$ based on
$( M^n, S, \nu)$. Simply, we call
the system (\ref{RSS}) a Liao perturbation of $(M^n, S, \nu)$.

From our Main Theorem we easily summarize the effect of Liao perturbations on
nonzero Lyapunov exponents
\begin{thm}
Let $S$ be a $C^1$ vector field on $M^n$ and let $\nu\in
E(M^n, \phi)$ be a probability that has $\ell$ nonzero Lyapunov exponents
$ \lambda_1 < \ldots < \lambda_\ell $
together with $n-\ell$ zero Lyapunov exponents.
Then there exists a $C^1$ neighborhood ${\mathcal X}^1$ of $S$
on the space of all $C^1$ vector fields on $M^n$, so that for each $X \in {\mathcal X}^1$,
its reduced standard system (\ref{RSS}) based on $(S, \nu)$
has $\lambda_1, \ldots , \lambda_\ell$ as Lyapunov exponents.
In other words, Liao perturbation preserves the nonzero Lyapunov exponents.
\end{thm}

\begin{rem}
Because the Lyapunov exponents are constant on $\Lambda$ and $F(\eta)\subset id_\ell(W)$,
from Proposition~\ref{prop:mu} and the Main Theorem, Theorem 6.1 is independent of the choice of $x\in \Lambda$ and
$\tilde \beta\in F(\eta)$ and thus the reduced standard system based on $(S, v)$.
\end{rem}

\baselineskip=0.45 true cm

\end{document}